# MINIMAL BANACH-TARSKI DECOMPOSITIONS

CESARE STRAFFELINI AND KILIAN ZAMBANINI

ABSTRACT. We investigate the problem of finding the minimum number of pieces necessary for dividing a three-dimensional sphere or a ball and reassembling it to form $n$ congruent copies of the original object, generalising a known result by Raphael Robinson.

## 1. INTRODUCTION

A century ago, Polish mathematicians Stefan Banach and Alfred Tarski managed to prove a result that, in some way, changed our view of mathematics forever. Their renowned theorem, published in [BT24], shows that it is possible to divide the closed unit ball $\mathbb{D}^3 \subseteq \mathbb{R}^3$ into finitely many pieces and reassemble them to form two balls of the same size, hence duplicating the total volume, something that in "real life" looks impossible.

Banach and Tarski's proof made use of the axiom of choice (AC), whose content is not constructive. The original goal of the two mathematicians was to provide supporting arguments to their aversion to AC, showing a paradoxical behaviour in its consequences. Indeed, there are alternative settings where mathematics can be developed, like the Solovay model introduced later in [Sol70], where the axiom of choice fails and every set of reals is Lebesgue measurable, hence the ball duplication result also fails.

Though, Banach and Tarski's result did not have the consequences they expected. Indeed, today the vast majority of the mathematical community accepts the axiom of choice and its consequences, considering them simply as counterintuitive but not as paradoxical.

Returning to [BT24], no estimate on the number of pieces necessary for the decomposition was given therein. In the following years many mathematicians tried to find an upper bound: John von Neumann claimed in [Neu29] that nine pieces suffice. Later, in [Sie45], Wacław Sierpiński proved that eight are enough. Eventually, in [Rob47], Raphael Robinson settled the question showing that five pieces are the minimal number.

In his proof, Robinson showed that the same construction for the unit sphere $\mathbb{S}^2$ (the boundary of $\mathbb{D}^3$) requires four pieces. Starting from such a minimal decomposition he was able to construct one for $\mathbb{D}^3$ with only five elements, one of which is a singleton, and he also showed that with less than five pieces it is impossible to reach the same result.

It is clear that, iterating this construction, one can start with a unit ball, divide it into finitely many pieces, and reassemble them into $n$ balls of the same size, for all natural numbers $n$. A question that now naturally arises is the following: what is the minimal quantity of pieces needed to produce $n$ unit balls starting with one?

This work is motivated by the authors' aim of answering the latter question. We started by looking whether it had already been discussed somewhere, and we only found a Mathematics StackExchange post from twelve years ago in which an anonymous user posed the same problem.

A tentative answer was suggested by Andrés Caicedo in [Cai13], where he claims that $5n-2$ pieces are surely sufficient and states that "it would be tempting to say that" the answer is $2n+1$, probably based on the idea that for the sphere $2n$ pieces are the correct value, and then one more piece is needed to pass from $\mathbb{S}^2$ to $\mathbb{D}^3$, as done by Robinson in the $n=2$ case.

In this work we show that, indeed, the necessary number of pieces for dividing $\mathbb{S}^2$ and reassembling it into $n$ copies is $2n$, while, for $\mathbb{D}^3$, we prove that the value of the optimal decomposition is $3n-1$ pieces, $n-1$ of which can be taken as singletons.

Concerning our contribution, we begin our proof with a group-theoretic result (Theorem 2.5), that builds upon the classical proof for the $n=2$ case (Theorem 2.2). In this regard see, for instance, the monograph [Wag85] or its updated edition [TW16].

We point out that each proof of Theorem 2.2 makes use of a suitable partition of $\mathbb{F}_2$ into four subsets. Each such partition is a slight variant of that formed by the sets "words beginning with $\sigma$", "words beginning with $\tau$", "words beginning with $\sigma^{-1}$" and "words beginning with $\tau^{-1}$", where $\sigma$ and $\tau$ are the generators of $\mathbb{F}_2$.

In the $n>2$ case, the main issue is to partition $\mathbb{F}_2$ into $2n$ subsets. As $\mathbb{F}_2$ is two-generated, finding a decomposition into $2n$ sets that is suitable for our later purposes is not completely trivial.

Another significant complication in the $n>2$ case arises in the proof of Theorem 4.4 below. Therein we need to single out a word $\omega$ in a very careful manner. Notice that, in the $n=2$ case, this is not needed because of the high symmetry of the four sets forming the $\mathbb{F}_2$ partition whose existence is claimed in Theorem 2.2. Such symmetry is then used throughout the whole proof.

## 2. Group Theory Prerequisites

For our results, we need basic information about the free group on two generators. We refer the reader to Chapter II, Section §5 of Paolo Aluffi's amazing algebra book [Alu09] for a detailed introduction and presentation of the topic. In the rest of this section we recall only the results needed in the following, in order to have a self-contained article.

**Definition 2.1.** Let $S$ be any set (of *letters*). We call *alphabet* on $S$ the set
$$S \sqcup \{\sigma^{-1} : \sigma \in S\},$$
where $\sigma^{-1}$ is just a formal notation, and $\sqcup$ denotes the disjoint union. Let $W$ the set of all finite juxtapositions of elements of the alphabet on $S$, namely,
$$W := \{\sigma_0^{e_0} \cdots \sigma_{m-1}^{e_{m-1}} : m \in \mathbb{N},\ \sigma_i \in S : i < m,\ e_i \in \{1,-1\} : i < m\},$$

where clearly $\sigma^1 := \sigma$ for $\sigma \in S$. The elements of $W$ are called *words* of the alphabet on $S$, and the set $W$ also contain the *empty word* (for $m = 0$), that we denote by $\varepsilon$.

We introduce an equivalence relation on $W$, saying that two words are equivalent if and only if we can reach one from the other by introducing and/or deleting a finite number of expressions of the form $\sigma\sigma^{-1}$ or $\sigma^{-1}\sigma$, for $\sigma \in S$. A word is said to be *reduced* if it contains no substring of the form $\sigma\sigma^{-1}$ or $\sigma^{-1}\sigma$, for any $\sigma \in S$.

Every word is equivalent to a unique reduced word, so the quotient of $W$ with respect to this equivalence relation can be identified with the set of reduced words. This set, denoted by $\mathbb{F}_S$, has a group structure with respect to the juxtaposition operation (up to equivalence), where the word $\varepsilon$ is the neutral element and
$$(\sigma_0^{e_0} \cdots \sigma_{m-1}^{e_{m-1}})^{-1} := \sigma_{m-1}^{-e_{m-1}} \cdots \sigma_0^{-e_0}.$$
This group is called the *free group* on the set $S$ of generators. If two sets $S$ and $S'$ are equipotent, then their corresponding free groups $\mathbb{F}_S$ and $\mathbb{F}_{S'}$ are isomorphic. Usually, $\mathbb{F}_k$ denotes "the" free group on $k$ generators, for some $k \in \mathbb{N}$.

From now onwards, we will only deal with the free group $\mathbb{F}_2$ on two generators, that we will call $\sigma$ and $\tau$. If $\gamma \in \mathbb{F}_2$ and $E \subseteq \mathbb{F}_2$, we let
$$\gamma(E) := \{\gamma \cdot \omega : \omega \in E\},$$
where $\gamma \cdot \omega$ (sometimes just $\gamma\omega$) denotes the juxtaposition of $\gamma$ and $\omega$.

In the classical proof of the Banach-Tarski decomposition theorem, with $n = 2$, the crucial group-theoretic result used therein is the following:

**Theorem 2.2.** *For each $\omega \in \mathbb{F}_2$, there is a partition of $\mathbb{F}_2$ in four subsets*
$$\mathbb{F}_2 = A_0 \sqcup B_0 \sqcup A_1 \sqcup B_1$$
*such that $\varepsilon$ and $\omega$ belong to the same element of the partition, and*
$$\sigma(B_0) = \mathbb{F}_2 - A_0; \quad \tau(B_1) = \mathbb{F}_2 - A_1.$$

We will not prove Theorem 2.2 because it follows from our Theorem 2.5 below. From now on, we fix a natural $n \geqslant 2$. This is the number of objects congruent to the original one we are building with our decompositions. If $n = 2$ we get a proof of the classical Banach-Tarski result.

Let $\alpha \in \mathbb{F}_2$ be a reduced word. We denote by $I(\alpha)$ the set of all reduced words in $\mathbb{F}_2$ that begin with $\alpha$. Then we define
$$A_0^* := I(\sigma), \ B_0^* := I(\sigma^{-n+1}), \ \gamma_0 := \sigma^{n-1};$$
$$A_1^* := I(\tau), \ B_1^* := I(\sigma^{-n+2}\tau^{-1}), \ \gamma_1 := \tau\sigma^{n-2};$$
and, for $2 \leqslant i < n$ (notice that in the $n = 2$ case there is no such $i$),
$$A_i^* := I(\sigma^{-i+1}\tau), \ B_i^* := I(\sigma^{-i+2}\tau^{-1}), \ \gamma_i := \sigma^{-i+1}\tau\sigma^{i-2}.$$

**Lemma 2.3.** *For all $i < n$, $\gamma_i(B_i^*) = \mathbb{F}_2 - A_i^*$ and $\gamma_i^{-1}(A_i^*) = \mathbb{F}_2 - B_i^*$.*

*Proof.* We start by showing that $\gamma_i(B_i^*) = \mathbb{F}_2 - A_i^*$ for all $i < n$. We distinguish the proof into three different cases: $i = 0$, $i = 1$ and $2 \leqslant i < n$.

If $i = 0$, let $\alpha \in B_0^*$, namely $\alpha = \sigma^{-n+1}\eta$, for some $\beta \in \mathbb{F}_2$ whose reduced word form does not begin with $\sigma$. Then
$$\gamma_0 \cdot \alpha = \sigma^{n-1}\sigma^{-n+1}\beta = \beta.$$

Since $\beta \notin I(\sigma)$, then $\gamma_0 \cdot \alpha = \beta \in \mathbb{F}_2 - A_0^*$, as desired. As for the converse inclusion, let $\alpha \in \mathbb{F}_2 - A_0^*$. This means that $\alpha$ does not begin with $\sigma$. We write
$$\alpha = \sigma^{n-1}\sigma^{-n+1}\alpha = \sigma^{n-1} \cdot \beta,$$
for some $\beta := \sigma^{-n+1}\alpha \in B_0^*$, hence $\alpha \in \gamma_0(B_0^*)$, as desired. The other cases, $i=1$ and $2 \leqslant i < n$, are similar, so we provide only brief sketches of their proofs.

If $\alpha \in B_1^*$, then $\alpha = \sigma^{-n+2}\tau^{-1}\beta$ for some $\beta \in \mathbb{F}_2$ not beginning with $\tau$. Then
$$\gamma_1 \cdot \alpha = \tau\sigma^{n-2}\sigma^{-n+2}\tau^{-1}\beta = \beta,$$
and $\beta \in \mathbb{F}_2 - A_1^*$. Conversely, if $\alpha \in \mathbb{F}_2 - A_1^*$, we can write
$$\alpha = \tau\sigma^{n-2}\sigma^{-n+2}\tau^{-1}\alpha = \tau\sigma^{n-2} \cdot \beta,$$
with $\beta := \sigma^{-n+2}\tau^{-1}\alpha$, and the latter is in $B_1^*$ since $\alpha \notin I(\tau)$. So, $\alpha \in \gamma_1(B_1^*)$.

Lastly, if $\alpha \in B_i^*$, for some $2 \leqslant i < n$, then $\alpha = \sigma^{-i+2}\tau^{-1}\beta$, where $\beta$ is a reduced word that does not begin with $\tau$. Hence it is straightforward to notice that
$$\gamma_i \cdot \alpha = \sigma^{-i+1}\tau\sigma^{i-2}\sigma^{-i+2}\tau^{-1}\beta = \sigma^{-i+1}\beta.$$
Since $\beta$ does not begin with $\tau$, we get that $\gamma_i \cdot \alpha \in \mathbb{F}_2 - A_i^*$, as desired. Regarding the reverse inclusion, let $\alpha \in \mathbb{F}_2 - A_i^*$. Then
$$\alpha = \sigma^{-i+1}\tau\sigma^{i-2}\sigma^{-i+2}\tau^{-1}\sigma^{i-1}\alpha = \sigma^{-i+1}\tau\sigma^{i-2} \cdot \beta,$$
with $\beta := \sigma^{-i+2}\tau^{-1}\sigma^{i-1}\alpha$, which belongs to $B_i^*$ if $\sigma^{i-1}\alpha$ does not start with $\tau$. But this is the case by the assumption $\alpha \in \mathbb{F}_2 - A_i^*$.

Since $\gamma_i(B_i^*) = \mathbb{F}_2 - A_i^*$, then it is straightforward that $\gamma_i^{-1}(A_i^*) = \mathbb{F}_2 - B_i^*$. □

**Lemma 2.4.** *Let $i < n$, $k \geqslant 1$ and $0 \leqslant m < n-1$. Then*
$$\gamma_i^k \sigma^{-m} \in A_i^* \text{ and } \gamma_i^{-k}\sigma^{-m} \in B_i^*.$$

*Proof.* Once again, we divide the proof into three cases: $i=0$, $i=1$ or $2 \leqslant i < n$.

Starting with the case $A_0^*$, we notice that
$$\gamma_0^k \sigma^{-m} = (\sigma^{n-1})^k \sigma^{-m} = \sigma^{nk-k-m}$$
and the facts that $m < n-1$ and $k \geqslant 1$ imply that the latter starts with $\sigma$, so we get that $\gamma_0^k\sigma^{-m} \in A_0^*$. As for $B_0^*$, we have that
$$\gamma_0^{-k}\sigma^{-m} = (\sigma^{n-1})^{-k}\sigma^{-m} = \sigma^{-nk+k-m}$$
begins with $\sigma^{-n+1}$, hence it belongs to $B_0^*$.

Concerning $A_1^*$, observe that
$$\gamma_1^k \sigma^{-m} = (\tau\sigma^{n-2})^k \sigma^{-m}$$
is in $A_1^*$. Regarding $B_1^*$, we notice that
$$\gamma_1^{-k}\sigma^{-m} = (\tau\sigma^{n-2})^{-k}\sigma^{-m} = (\sigma^{-n+2}\tau^{-1})^k\sigma^{-m},$$
which starts by $\sigma^{-n+2}\tau^{-1}$ and belongs to $B_1^*$.

If instead $2 \leqslant i < n$, we get that
$$\gamma_i^k \sigma^{-m} = (\sigma^{-i+1}\tau\sigma^{i-2})^k \sigma^{-m},$$
which starts with $\sigma^{-i+1}\tau$ and therefore lies in $A_i^*$. To conclude, if $2 \leqslant i < n$, then
$$\gamma_i^{-k}\sigma^{-m} = (\sigma^{-i+1}\tau\sigma^{i-2})^{-k}\sigma^{-m} = (\sigma^{-i+2}\tau^{-1}\sigma^{i-1})^k\sigma^{-m}.$$

Such a word begins with $\sigma^{-i+2}\tau^{-1}$, therefore it is in $B_i^*$. □

We are ready to prove Theorem 2.5, generalising the classical Theorem 2.2.

Notice that in the case $n=2$ the coefficients $\gamma_0$ and $\gamma_1$ are just $\sigma$ and $\tau$. Hence the proof of Theorem 2.2 is easier because it is more "symmetric", in the sense that there are simply definable automorphisms of $\mathbb{F}_2$ with $\sigma \mapsto \tau$, with $\sigma \mapsto \sigma^{-1}$ or with $\sigma \mapsto \tau^{-1}$ that preserve the family $\{A_0^*, B_0^*, A_1^*, B_1^*\} = \{I(\sigma), I(\tau), I(\sigma^{-1}), I(\tau^{-1})\}$.

By this symmetric property, in the proof of Theorem 2.2 it is enough to assume, for example, $\omega \in I(\sigma)$. Here instead we consider all cases of $\omega \in A_i^*$ or $\omega \in B_i^*$.

**Theorem 2.5.** *For each $\omega \in \mathbb{F}_2$, there is a partition of $\mathbb{F}_2$ in $2n$-many subsets*

$$\mathbb{F}_2 = A_0 \sqcup B_0 \sqcup \ldots \sqcup A_{n-1} \sqcup B_{n-1}$$

*with $\varepsilon$ and $\omega$ in the same element of the partition, and such that for all $i<n$*

$$\gamma_i(B_i) = \mathbb{F}_2 - A_i.$$

*Proof.* First of all, we notice that, with the notation of Lemma 2.3 in force, we have

$$A_0^* \sqcup B_0^* \sqcup \ldots \sqcup A_{n-1}^* \sqcup B_{n-1}^* = \mathbb{F}_2 - \{\varepsilon, \sigma^{-1}, \ldots, \sigma^{-n+2}\}.$$

If $\omega \in \{\varepsilon, \sigma^{-1}, \ldots, \sigma^{-n+2}\}$, we will just behave as in the case $\omega = \sigma^{-n+1}$.

Let $\omega \notin \{\varepsilon, \sigma^{-1}, \ldots, \sigma^{-n+2}\}$. Therefore $\omega$ belongs to $A_i^*$ or $B_i^*$ for exactly one $i < n$. Let $C^*$ denote this set, hence $C^* := A_i^*$ if $\omega \in A_i^*$ and $C^* := B_i^*$ if $\omega \in B_i^*$.

Denote by $D^*$ the other one, so $D^* := B_i^*$ if $C^* = A_i^*$ and $D^* := A_i^*$ if $C^* = B_i^*$. Let also $\delta \in \mathbb{F}_2$ denote the element of the form $\gamma_i$ or $\gamma_i^{-1}$ such that $\delta(D^*) = \mathbb{F}_2 - C^*$, whose existence is ensured by Lemma 2.3.

We define $A_j := A_j^*$ and $B_j := B_j^*$ for $j \neq i$, while

$$C := C^* \sqcup \{\delta^{-k}\sigma^{-m} : 0 \leqslant m < n-1, k \geqslant 0\},$$
$$D := D^* - \{\delta^{-k}\sigma^{-m} : 0 \leqslant m < n-1, k \geqslant 1\}.$$

Finally, let $A_i := C$, $B_i := D$ if $\omega \in A_i^*$ and $A_i := D$, $B_i := C$ if $\omega \in B_i^*$.

Notice that, by Lemma 2.4, $\{\delta^{-k}\sigma^{-m} : 0 \leqslant m < n-1, k \geqslant 1\} \subseteq D^*$. Hence,

$$A_0 \sqcup B_0 \sqcup \ldots \sqcup A_{n-1} \sqcup B_{n-1} = \mathbb{F}_2$$

and $\omega$ and $\varepsilon$ are in the same element of the partition. From Lemma 2.3 it follows that $\gamma_j(B_j) = \mathbb{F}_2 - A_j$ for $j \neq i$. So it remains to prove that $\delta(D) = \mathbb{F}_2 - C$.

Since we know that $\delta(D^*) = \mathbb{F}_2 - C^*$, our desired conclusion follows by definition of $C$ and $D$, from Lemma 2.4, and by the trivial observation that

$$\delta(\{\delta^{-k}\sigma^{-m} : 0 \leqslant m < n-1, k \geqslant 1\}) = \{\delta^{-k}\sigma^{-m} : 0 \leqslant m < n-1, k \geqslant 0\}. \quad \square$$

We anticipate that in the proof of our generalisation we also need a refinement of Theorem 2.5, in which we slightly change the definition of the subsets $A_j$ in the case where the word $\omega$ begins with the letter $\tau$, again showcasing the asymmetry between the two generators $\sigma$ and $\tau$ that was absent in the $n=2$ case.

**Lemma 2.6.** *If $\omega \in \mathbb{F}_2$ begins with $\tau$, there is a partition $\langle A_i, B_i : 0 \leqslant i < n \rangle$ of $\mathbb{F}_2$ in $2n$-many sets such that $\sigma^{-i+1} \in A_i$ for all $i < n$, $\omega \in A_1$ and, for all $i < n$,*

$$\gamma_i(B_i) = \mathbb{F}_2 - A_i.$$

*Proof.* With the notation of Lemma 2.3 in force, same as in Theorem 2.5 we have

$$A_0^* \sqcup B_0^* \sqcup \ldots \sqcup A_{n-1}^* \sqcup B_{n-1}^* = \mathbb{F}_2 - \{\varepsilon, \sigma^{-1}, \ldots, \sigma^{-n+2}\}.$$

We let $A_0 := A_0^*$ (notice that $\sigma^1 = \sigma \in A_0$, as desired), $B_0 := B_0^*$;

$$A_1 := A_1^* \sqcup \{\gamma_1^{-k} : k \geqslant 0\}, \quad B_1 := B_1^* - \{\gamma_1^{-k} : k \geqslant 1\}.$$

We also notice that $\sigma^0 = \varepsilon = \gamma_1^{-0} \in A_1$, as desired.

Also, by Lemma 2.4, we have

$$A_1 \sqcup B_1 = A_1^* \sqcup B_1^* \sqcup \{\varepsilon\}.$$

Finally, by Lemma 2.3 and by the trivial observation that

$$\gamma_1(\{\gamma_1^{-k} : k \geqslant 1\}) = \{\gamma_1^{-k} : k \geqslant 0\}$$

we get that $\gamma_1(B_1) = \mathbb{F}_2 - A_1$.

In a very similar way, for all $2 \leqslant i < n$ we define

$$A_i := A_i^* \sqcup \{\gamma_i^{-k}\sigma^{-i+1} : k \geqslant 0\}, \quad B_i := B_i^* - \{\gamma_i^{-k}\sigma^{-i+1} : k \geqslant 1\}.$$

Thus $\sigma^{-i+1} = \gamma_i^0 \sigma^{-i+1} \in A_i$, as required, and by Lemma 2.4 we have

$$A_i \sqcup B_i = A_i^* \sqcup B_i^* \sqcup \{\sigma^{-i+1}\}.$$

Finally by Lemma 2.3 and by the equality

$$\gamma_i(\{\gamma_i^{-k}\sigma^{-i+1} : k \geqslant 1\}) = \{\gamma_i^{-k}\sigma^{-i+1} : k \geqslant 0\}$$

we get $\gamma_i(B_i) = \mathbb{F}_2 - A_i$. Putting everything together, we conclude. □

## 3. Paradoxality under Group Actions

In the next results, we deal with group actions. We refer to Chapter II, Section §9 of [Alu09] for a thorough presentation. We just recall some of the basic definitions.

If there is no ambiguity, we identify a group $\langle G, \cdot, \varepsilon \rangle$ with its support $G$.

**Definition 3.1.** Let $G$ be a group and $\Omega$ be a nonempty set. A *group action* of $G$ on $\Omega$ is a map from $G \times \Omega$ to $\Omega$, denoted by $\langle \gamma, x \rangle \mapsto \gamma \bullet x$, such that

$$\varepsilon \bullet x = x, \text{ and } (\gamma \cdot \delta) \bullet x = \gamma \bullet (\delta \bullet x)$$

for all $\gamma$ and $\delta$ in $G$. Given a group action of $G$ on $\Omega$, we say that $G$ *acts* on $\Omega$ and we write $G \curvearrowright \Omega$. Every group $G$ acts on its support $G$ with the action given by

$$\gamma \bullet \delta := \gamma \cdot \delta.$$

We call this the *natural* action of $G$ on itself, and we denote it as $G \curvearrowright G$.

If $G \curvearrowright \Omega$ is a group action, $\gamma \in G$ and $E \subseteq \Omega$, we let

$$\gamma(E) := \{\gamma \bullet x : x \in E\}$$

The latter is consistent with the notation introduced after Definition 2.1 when we consider the natural action $G \curvearrowright G$.

**Definition 3.2.** Let $G \curvearrowright \Omega$, $E$ and $F$ subsets of $\Omega$ and $r \in \mathbb{N}$. We say that $E$ and $F$ are *G-equidecomposable* with $r$ pieces, and we write
$$E \sim_r F,$$
if there are a partition
$$E = E_0 \sqcup \ldots \sqcup E_{r-1}$$
and a sequence $\langle \gamma_j : j < r \rangle$ of elements of $G$ with
$$F = \gamma_0(E_0) \sqcup \ldots \sqcup \gamma_{r-1}(E_{r-1}).$$

**Definition 3.3.** Let $G \curvearrowright \Omega$ and let $E \subseteq \Omega$. We say $E$ is $(n,r)$-*paradoxical* with respect to the action $G \curvearrowright \Omega$ if there is a partition $E = P_0 \sqcup \ldots \sqcup P_{n-1}$ such that
$$P_j \sim_{r_j} E \text{ for all } j < n, \text{ and } r_0 + \ldots + r_{n-1} = r.$$

If the support of a group $G$ is $(n,r)$-paradoxical with respect to its own natural action $G \curvearrowright G$, then we will just say that the group $G$ is $(n,r)$-paradoxical.

**Theorem 3.4.** *If $\Omega$ is $(n,r)$-paradoxical with respect to $G \curvearrowright \Omega$, then $r \geqslant 2n$.*

*Proof.* For sake of contradiction, we assume that $r < 2n$. Let $\langle P_j : j < n \rangle$ be a partition of $\Omega$ and $\langle r_j : j < n \rangle$ be a sequence of positive natural numbers such that
$$P_j \sim_{r_j} \Omega \text{ for all } j < n, \text{ and } r_0 + \ldots + r_{n-1} = r.$$
Since $r < 2n$ and $r_j \geqslant 1$ for all $j < n$, there is at least one $j_0 < n$ such that $r_{j_0} = 1$.

Hence $P_{j_0} \sim_1 \Omega$. Namely, there is some $\gamma_0 \in G$ such that
$$\Omega = \gamma_0(P_{j_0}).$$
This implies $\gamma_0^{-1}(\Omega) = P_{j_0}$, but we also know that $\gamma_0^{-1}(\Omega) = \Omega$. Hence $P_j = \varnothing$ for all $j \neq j_0$, a contradiction with $P_j \sim_{r_j} \Omega$ by the assumption that $n \geqslant 2$. □

As a corollary to the latter and to Theorem 2.5, we now get the following result:

**Corollary 3.5.** $\mathbb{F}_2$ *is $(n, 2n)$-paradoxical and not $(n,r)$-paradoxical for any $r < 2n$.*

*Proof.* By Theorem 3.4, we know that $\mathbb{F}_2$ is not $(n,r)$-paradoxical for $r < 2n$. Using the notation of Theorem 2.5, let $P_j := A_j \sqcup B_j$ for $j < n$. Then $\langle P_j : j < n \rangle$ is a partition of $\mathbb{F}_2$. Moreover, $P_j \sim_2 \mathbb{F}_2$ for all $j < n$, since
$$\mathbb{F}_2 = A_j \sqcup (\mathbb{F}_2 - A_j) = \varepsilon(A_j) \sqcup \gamma_j(B_j). \qquad \square$$

Now we are ready to start proving paradoxality results for $\mathbb{S}^2$ and $\mathbb{D}^3$. The idea is to *transfer* the paradoxical behaviour of $\mathbb{F}_2$ to other sets on which $\mathbb{F}_2$ acts.

**Definition 3.6.** Let $G \curvearrowright \Omega$ and let $x \in \Omega$. The *stabiliser* of $x$ is
$$\mathrm{Stab}(x) := \{\gamma \in G : \gamma \bullet x = x\}.$$
We say that $x$ has *trivial* stabiliser if $\mathrm{Stab}(x) = \{\varepsilon\}$. Instead, the *orbit* of $x$ is
$$\mathcal{O}(x) := \{\delta \bullet x : \delta \in G\}.$$
Notice that if $y \in \mathcal{O}(x)$, namely, $y = \delta \bullet x$ for some $\delta \in G$, then
$$\mathrm{Stab}(y) = \{\gamma \in G : \delta^{-1}\gamma\delta \in \mathrm{Stab}(x)\}.$$
In particular, the property of "having a trivial stabiliser" is constant on orbits.

An action $G \curvearrowright \Omega$ is said to be *free* if $\mathrm{Stab}(x)$ is trivial for all $x \in \Omega$.

**Theorem 3.7.** *Let $G \curvearrowright \Omega$, with a free action. Suppose that $G$ is $(n,r)$-paradoxical. Then $\Omega$ is $(n,r)$-paradoxical with respect to the action $G \curvearrowright \Omega$.*

*Proof.* Using AC, we pick an element from each orbit of the action of $G$ on $\Omega$ and collect them into a set $M$. If $S \subseteq G$, we let $S' := \{\gamma \bullet m : \gamma \in S, m \in M\}$.

Since $G$ is $(n,r)$-paradoxical, there exist subsets $\langle P_j : j < n \rangle$ such that
$$G = P_0 \sqcup \ldots \sqcup P_{n-1}$$
and with the property that
$$P_j \sim_{r_j} G \text{ for } j < n \text{ and } r_0 + \ldots + r_{n-1} = r.$$
We claim that $\Omega = P'_0 \sqcup \ldots \sqcup P'_{n-1}$ and that $P'_j \sim_{r_j} \Omega$ for all $j < n$.

First of all, notice that $\langle P'_j : j < n \rangle$ is a partition of $\Omega$. Indeed, if $y \in \Omega$, then $y = \gamma \bullet m$, for some $\gamma \in G$, $m \in M$. Since $G = P_0 \cup \ldots \cup P_{n-1}$, we get that
$$\Omega = P'_0 \cup \ldots \cup P'_{n-1}.$$
Moreover, if $y \in P'_i \cap P'_j$, then there exist $\gamma_i \in P_i$, $\gamma_j \in P_j$, $m_i \in M$, $m_j \in M$ with $\gamma_i \bullet m_i = y = \gamma_j \bullet m_j$. This implies $(\gamma_i^{-1}\gamma_j) \bullet m_j = m_i$ and, by definition of $M$, $m_i = m_j$. Since the action of $G$ on $\Omega$ is free, we conclude that $\gamma_i = \gamma_j$. Since $\langle P_j : j < n \rangle$ is a partition of $G$, then $i = j$. Hence $P'_i \cap P'_j = \varnothing$ if $i \neq j$.

Next we fix $j < n$: our goal is to prove that $P'_j \sim_{r_j} \Omega$. Since $P_j \sim_{r_j} G$, there exists a partition $P_j = E_0 \sqcup \ldots \sqcup E_{r_j - 1}$ and some $\langle \gamma_k : k < r_j \rangle$ in $G$ such that
$$G = \gamma_0(E_0) \sqcup \ldots \sqcup \gamma_{r_j - 1}(E_{r_j - 1}).$$
By repeating the argument above, we obtain that $P'_j = E'_0 \sqcup \ldots \sqcup E'_{r_j - 1}$.

It remains to show that $\Omega = \gamma_0(E'_0) \sqcup \ldots \sqcup \gamma_{r_j - 1}(E'_{r_j - 1})$. To do so, we fix $y \in \Omega$. Then $y = \gamma \bullet m$ for some $\gamma \in G$ and $m \in M$. But $G = \gamma_0(E_0) \sqcup \ldots \sqcup \gamma_{r_j - 1}(E_{r_j - 1})$, so $y = (\gamma_i \delta_i) \bullet m = \gamma_i \bullet (\delta_i \bullet m)$ for some $i < r_j$ and $\delta_i \in E_i$. To conclude, we notice that if $y \in \gamma_i(E'_i) \cap \gamma_j(E'_j)$ there are $\delta_i \in E_i$, $\delta_j \in E_j$, $m_i \in M$, $m_j \in M$ with
$$\gamma_i \bullet (\delta_i \bullet m_i) = y = \gamma_j \bullet (\delta_j \bullet m_j).$$
By the axioms of group action and by the triviality of the stabilisers, it follows that $\gamma_i \delta_i = \gamma_j \delta_j$, implying $i = j$ since $G = \gamma_0(E_0) \sqcup \ldots \sqcup \gamma_{r_j - 1}(E_{r_j - 1})$. □

Note that in the proof of Theorem 3.7, as in the proof of Theorem 4.4 below, we make use of AC to pick one element from each orbit of the action $G \curvearrowright \Omega$.

However, exactly as proven in Janusz Pawlikowski's article [Paw91], it is enough to assume the Hahn-Banach theorem, a strictly weaker version of AC.

As we already mentioned in the introduction, it is impossible to prove this kind of results without appealing to any fragment of the axiom of choice.

We recall a well-known notation in group theory: given two groups $G$ and $H$, we write $H \leqslant G$ if $H$ is a subgroup of $G$ and $G \simeq H$ if $G$ and $H$ are isomorphic.

**Corollary 3.8.** *If $H \leqslant G$ is such that $H \simeq \mathbb{F}_2$, then $G$ is $(n, 2n)$-paradoxical.*

*Proof.* By Corollary 3.5, $H$ is $(n, 2n)$-paradoxical. Since the action $H \curvearrowright G$ by left multiplication is free, by Theorem 3.7 $G$ is $(n, 2n)$-paradoxical with respect to this action. Since $H \leqslant G$, we deduce that $G$ is itself $(n, 2n)$-paradoxical. □

## 4. The Result for the Sphere

In this section we would like to apply the results from Section §3 to the case of the special orthogonal group $\mathrm{SO}(3)$ acting on the sphere $\mathbb{S}^2 \subseteq \mathbb{R}^3$.

We quote the following, a proof of which can be found in [Świ94]:

**Fact 4.1.** *The group* $\mathrm{SO}(3)$ *admits a subgroup isomorphic to* $\mathbb{F}_2$.

If the action of $\mathrm{SO}(3)$ on $\mathbb{S}^2$ were free, by the previous results we would get the paradoxality of the sphere. Unfortunately, every rotation in $\mathrm{SO}(3)$ different from the identity has two fixed points on $\mathbb{S}^2$, hence there are nontrivial stabilisers.

Nevertheless, such an action satisfies the following property:

**Definition 4.2.** A group action $G \curvearrowright \Omega$ is *locally abelian* if the subgroup $\mathrm{Stab}(x) \leqslant G$ is abelian for all elements $x \in \Omega$.

**Lemma 4.3.** *The action of* $\mathrm{SO}(3)$ *on* $\mathbb{S}^2$ *is locally abelian.*

*Proof.* If two elements of $\mathrm{SO}(3)$ fix a common point in $\mathbb{R}^3$ which is different from the origin, then they have the same rotation axis, hence they commute. □

The following result is crucial:

**Theorem 4.4.** *If the action* $\mathbb{F}_2 \curvearrowright \Omega$ *is locally abelian, then* $\Omega$ *is* $(n, 2n)$-*paradoxical with respect to this action.*

*Proof.* Let $\sigma$ and $\tau$ be generators of $\mathbb{F}_2$. We prove the existence of a partition of $\Omega$ witnessing its $(n, 2n)$-paradoxicality. Since the orbits of the action $\mathbb{F}_2 \curvearrowright \Omega$ form a partition of $\Omega$, we work on each orbit. Let $\gamma_i \in \mathbb{F}_2$ be defined, as in Section §2, as

$$\gamma_0 := \sigma^{n-1}; \quad \gamma_1 := \tau\sigma^{n-2}; \quad \gamma_i := \sigma^{-i+1}\tau\sigma^{i-2} \text{ for } 2 \leqslant i < n.$$

For each orbit $\mathcal{O}$, we prove that we can write

$$\mathcal{O} = A_0^{\mathcal{O}} \sqcup \cdots \sqcup A_{n-1}^{\mathcal{O}} \sqcup B_0^{\mathcal{O}} \sqcup \cdots \sqcup B_{n-1}^{\mathcal{O}}$$

for some $A_i^{\mathcal{O}}$ and $B_i^{\mathcal{O}}$ for $i < n$ with the property $\gamma_i(B_i^{\mathcal{O}}) = \mathcal{O} - A_i^{\mathcal{O}}$ for every $i < n$. We then let $A_i' := \bigcup_{\mathcal{O}} A_i^{\mathcal{O}}$, $B_i' := \bigcup_{\mathcal{O}} B_i^{\mathcal{O}}$ to get the conclusion for the entire set $\Omega$.

As remarked before, given an orbit $\mathcal{O}$, either the stabilisers of elements of $\mathcal{O}$ are all trivial or they are all nontrivial. First we consider the case of an orbit $\mathcal{O}$ whose elements have trivial stabilisers. We pick $x \in \mathcal{O}$. Since the stabiliser of $x$ is trivial, each element $y \in \mathcal{O}$ can be written in a *unique* way as $y = \alpha \bullet x$, for some $\alpha \in \mathbb{F}_2$.

By Theorem 2.5, we can write $\mathbb{F}_2 = A_0 \sqcup B_0 \sqcup \ldots \sqcup A_{n-1} \sqcup B_{n-1}$, with

$$\gamma_i(B_i) = \mathbb{F}_2 - A_i$$

for all $i < n$. We define $A_i^{\mathcal{O}}$ and $B_i^{\mathcal{O}}$ in the following way: we stipulate that $y \in A_i^{\mathcal{O}}$ if and only if $\alpha \in A_i$, $y \in B_i^{\mathcal{O}}$ if and only if $\alpha \in B_i$.

By construction, it is clear that

$$\mathcal{O} = A_0^{\mathcal{O}} \sqcup \cdots \sqcup A_{n-1}^{\mathcal{O}} \sqcup B_0^{\mathcal{O}} \sqcup \cdots \sqcup B_{n-1}^{\mathcal{O}}.$$

Moreover, if $z = \gamma_i \bullet y$ for some $y \in B_i^{\mathcal{O}}$, then $z = \gamma_i \bullet (\alpha \bullet x) = (\gamma_i\alpha) \bullet x$ for some $\alpha \in B_i$, which implies $z \in \mathcal{O} - A_i^{\mathcal{O}}$ by the property $\gamma_i(B_i) = \mathbb{F}_2 - A_i$. Conversely, if $z \in \mathcal{O} - A_i^{\mathcal{O}}$, then $z = \alpha \bullet x$ for some $\alpha \in \mathbb{F}_2 - A_i = \gamma_i(B_i)$. Hence $z = \gamma_i \bullet (\beta \bullet x)$ for some $\beta \in B_i$. Hence $z \in \gamma_i(B_i^{\mathcal{O}})$, ending the proof of the first case.

Next, let $\mathcal{O}$ be an orbit whose elements have nontrivial stabilisers. Let $\omega \in \mathbb{F}_2$, $\omega \neq \varepsilon$, be of minimal length among the words in $\mathbb{F}_2$ that fix at least one element of $\mathcal{O}$. We also pick an element $x \in \mathcal{O}$ such that $\omega \bullet x = x$. Let $\alpha$ be the initial letter of the word $\omega$. Notice that $\omega$ does not end with $\alpha^{-1}$, otherwise from

$$(\alpha^{-1}\omega\alpha) \bullet (\alpha^{-1} \bullet x) = (\alpha^{-1}\omega) \bullet x = \alpha^{-1} \bullet x,$$

we would get that the word $\alpha^{-1}\omega\alpha$, which has length strictly smaller than the length of $\omega$, fixes the element $\alpha^{-1} \bullet x$, violating the minimality of $\omega$.

Notice that if $\omega \bullet x = x$ and $\alpha$ is as above, then $\omega' := \alpha^{-1}\omega\alpha$ has the same length as $\omega$ and fixes $\alpha^{-1} \bullet x$, as we have already pointed out. A similar argument holds for $\omega^{-1}$: it has the same length as $\omega$ and fixes $x$. Therefore, if $\omega$ contains the letter $\tau$, up to replacing $\omega$ with $\omega'$ and possibly repeating the previous argument, we can assume that $\omega$ begins with $\tau$ (so, $\alpha = \tau$). Similarly, if $\omega$ contains the letter $\tau^{-1}$, then we can replace $\omega$ with $\omega^{-1}$ and then make use of the previous argument.

Eventually, if $\omega$ does not contain neither $\tau$ nor $\tau^{-1}$, then it is a power of $\sigma$. Up to replacing $\omega$ with $\omega^{-1}$, we can assume that $\omega = \sigma^{-k}$ for some $k > 0$.

To finish the proof, we distinguish in the two cases whether $\alpha = \tau$ or $\omega = \sigma^{-k}$. This is the content of Lemmata 4.6 and 4.7 below, which we prove next. □

Before proving Lemmata 4.6 and 4.7, we show the preliminary result Lemma 4.5 concerning locally abelian actions. The proof is the same as that of the case $n = 2$. See, for instance, [Wag85]. We present it in detail, for sake of completeness.

**Lemma 4.5.** *Let $\mathbb{F}_2 \curvearrowright \Omega$ be locally abelian, $\mathcal{O}$ an orbit with nontrivial stabilisers, and let $\omega \neq \varepsilon$ be of minimal length fixing at least one element $x \in \mathcal{O}$. Let $\alpha$ be the first letter of $\omega$. Then, for each $y \in \mathcal{O}$, there is a unique reduced word $\zeta \in \mathbb{F}_2$ that does not end either in $\alpha^{-1}$ or in $\omega$ and such that $y = \zeta \bullet x$.*

*Proof.* We first prove existence and then uniqueness.

EXISTENCE. Let $y \in \mathcal{O}$. Then there is some $\varrho \in \mathbb{F}_2$ such that $y = \varrho \bullet x$. If $\varrho$ does not end in $\omega$ or in $\alpha^{-1}$, we are done. If $\varrho$ ends in $\omega$, we consider $\varrho\omega^{-1}$, which still has the property that $y = (\varrho\omega^{-1}) \bullet x$, since $\omega$ fixes $x$, and we can repeat the process until $\varrho\omega^{-k}$ does not end in $\omega$. Notice that in this case $\varrho\omega^{-k}$ does not end in $\alpha^{-1}$, because otherwise $\varrho$ would have a substring $\alpha^{-1}\alpha$ in its reduced form.

If $\varrho$ ends in $\alpha^{-1}$, we consider $\varrho\omega$. Hence $y = (\varrho\omega) \bullet x$, and we can repeat the argument until $\varrho\omega^k$ does not end in $\alpha^{-1}$ (recall that $\omega$ does not end in $\alpha^{-1}$). Notice also that, by construction, this new word does not end in $\omega$.

UNIQUENESS. First we prove that the only elements of $\mathbb{F}_2$ that fix $x$ are the powers of $\omega$. Indeed, if $\psi$ also fixes $x$, then $\psi\omega = \omega\psi$ by the local abelian property of the action of $\mathbb{F}_2$ on $\Omega$. However, if two elements of a free group commute, then they are powers of a common element. This means that there exist $\phi \in \mathbb{F}_2$ and $j, k \in \mathbb{Z}$ such that $\omega = \phi^j$, $\psi = \phi^k$. By minimality of length of $\omega$ it follows $|k| \geq |j|$. Let $k = j\ell + r$, with $\ell \in \mathbb{Z}$ and $0 \leq r < |j|$. Then

$$x = \psi \bullet x = \phi^k \bullet x = \phi^r(\phi^j)^\ell \bullet x = \phi^r\omega^\ell \bullet x = \phi^r \bullet x.$$

So $\phi^r$ fixes $x$ and, by minimality of $\omega$, we get $r = 0$, forcing $\psi$ to be a power of $\omega$.

Let now assume that $y = \eta \bullet x = \zeta \bullet x$, for some $\eta, \zeta \in \mathbb{F}_2$ not ending neither in $\alpha^{-1}$ nor in $\omega$. Then $(\eta^{-1}\zeta) \bullet x = x$, so $\eta^{-1}\zeta$ is a power of $\omega$.

If $\eta^{-1}\zeta$ is a positive power of $\omega$ and $\eta^{-1}$ begins with $\alpha$, then $\eta$ ends with $\alpha^{-1}$: a contradiction. If instead $\eta^{-1}\zeta$ is a positive power of $\omega$ and $\eta^{-1}$ does not start with $\alpha$, notice that $\eta^{-1}\zeta$ starts with $\alpha$. Therefore $\eta^{-1}$ cancels with the first letters of $\zeta$ and so $\zeta$ does end with $\omega$: again, a contradiction. If $\eta^{-1}\zeta$ is a negative power of $\omega$, $\zeta^{-1}\eta$ is a positive power of $\omega$ and we reach a contradiction by the same argument.

Hence the only possibility is that $\eta^{-1}\zeta = \varepsilon$. So uniqueness is proved. $\square$

**Lemma 4.6.** *Let* $\mathbb{F}_2 \curvearrowright \Omega$ *be locally abelian,* $\mathcal{O}$ *an orbit with nontrivial stabilisers, and let* $\omega \neq \varepsilon$ *be of minimal length fixing at least one element* $x \in \mathcal{O}$.

*Assume that* $\omega$ *begins with* $\tau$. *Then there is a partition*
$$\mathcal{O} = A_0^{\mathcal{O}} \sqcup \cdots \sqcup A_{n-1}^{\mathcal{O}} \sqcup B_0^{\mathcal{O}} \sqcup \cdots \sqcup B_{n-1}^{\mathcal{O}}$$
*for some* $A_i^{\mathcal{O}}$ *and* $B_i^{\mathcal{O}}$ *for* $i < n$ *such that* $\gamma_i(B_i^{\mathcal{O}}) = \mathcal{O} - A_i^{\mathcal{O}}$ *for every* $i < n$.

*Proof.* We apply the construction of Lemma 2.6 to get a partition
$$\mathbb{F}_2 = A_0 \sqcup B_0 \sqcup \ldots \sqcup A_{n-1} \sqcup B_{n-1}$$
such that $\sigma^{-i+1} \in A_i$, for all $i < n$, $\omega \in A_1$ and, for all $i < n$,
$$\gamma_i(B_i) = \mathbb{F}_2 - A_i.$$
Next we define the sets $A_i^{\mathcal{O}}$ and $B_i^{\mathcal{O}}$, for $0 \leqslant i < n$, in the following way: for each $y \in \mathcal{O}$, let $\zeta$ be the word as in Lemma 4.5. We stipulate that $y \in A_i^{\mathcal{O}}$ or $y \in B_i^{\mathcal{O}}$ if and only if $\zeta \in A_i$ or $\zeta \in B_i$, respectively. Since
$$\mathbb{F}_2 = A_0 \sqcup \ldots \sqcup A_{n-1} \sqcup B_0 \sqcup \ldots \sqcup B_{n-1},$$
by uniqueness of $\zeta$,
$$\mathcal{O} = A_0^{\mathcal{O}} \sqcup \ldots \sqcup A_{n-1}^{\mathcal{O}} \sqcup B_0^{\mathcal{O}} \sqcup \ldots \sqcup B_{n-1}^{\mathcal{O}}.$$
It remains to prove that $\gamma_i(B_i^{\mathcal{O}}) = \mathcal{O} - A_i^{\mathcal{O}}$, for all $i < n$.

First we prove that $\gamma_i(B_i^{\mathcal{O}}) \subseteq \mathcal{O} - A_i^{\mathcal{O}}$, for all $i < n$. Let $y \in B_i^{\mathcal{O}}$, hence $y = \zeta \bullet x$ for some $\zeta \in B_i$ that does not end neither in $\tau^{-1}$ nor in $\omega$. Since
$$\gamma_i(B_i) = \mathbb{F}_2 - A_i,$$
we can write $\gamma_i \bullet y = (\gamma_i \zeta) \bullet x$, with $\gamma_i \zeta \in \mathbb{F}_2 - A_i$. In particular, if $\gamma_i \zeta$ does not end neither in $\tau^{-1}$ nor in $\omega$, we conclude that $\gamma_i \bullet y \in \mathcal{O} - A_i^{\mathcal{O}}$.

If $i = 0$, then $\gamma_0 = \sigma^{n-1}$, which does not contain neither $\tau$ nor $\tau^{-1}$. Since we are assuming that $\zeta$ does not end neither in $\tau^{-1}$ nor in $\omega$, we deduce that the same holds for $\gamma_0 \zeta = \sigma^{n-1}\zeta$ (recall that $\omega$ begins with $\tau$), as desired.

If $i = 1$, then $\gamma_1 = \tau\sigma^{n-2}$. Since $\zeta$ does not end in $\tau^{-1}$, the same holds for
$$\gamma_1 \zeta = \tau\sigma^{n-2}\zeta.$$
Moreover, since $\zeta$ does not end in $\omega$ and $\omega$ begins with $\tau$, the only possible case in which $\tau\sigma^{n-2}\zeta$ ends in $\omega$ is when $\tau\sigma^{n-2}\zeta = \omega$. In such case $\omega = \gamma_1\zeta \in \mathbb{F}_2 - A_1$, which contradicts $\omega \in A_1$ as follows from the construction in Lemma 2.6.

If $2 \leqslant i < n$, then $\gamma_i = \sigma^{-i+1}\tau\sigma^{i-2}$ and, as before, it is easy to see that $\gamma_i \zeta$ does not end in $\tau^{-1}$ if $\zeta$ does not. Finally, we have to exclude the case in which
$$\gamma_i \zeta = \sigma^{-i+1}\tau\sigma^{i-2}\zeta$$
ends in $\omega$. Since $\zeta$ does not end in $\omega$, and $\omega$ begins with $\tau$, then the only possibility would be $\omega = \tau\sigma^{i-2}\zeta$. However, in this case, $\gamma_i\zeta = \sigma^{-i+1}\omega$ begins with $\sigma^{-i+1}\tau$, so it belongs to $A_i^* \subseteq A_i$ (for the definition of $A_i^*$, refer to the notation of Lemma 2.3 and the construction in Lemma 2.6), contradicting the property that $\gamma_i \zeta \in \mathbb{F}_2 - A_i$.

Eventually, we are left to prove $\gamma_i(B_i^{\mathcal{O}}) \supseteq \mathcal{O} - A_i^{\mathcal{O}}$ for $i < n$, equivalent to
$$\gamma_i^{-1}(\mathcal{O} - A_i^{\mathcal{O}}) \subseteq B_i^{\mathcal{O}}.$$
Let $y \in \mathcal{O} - A_i^{\mathcal{O}}$. Hence $y = \zeta \bullet x$ for some $\zeta \in \mathbb{F}_2 - A_i$ which does not end neither in $\tau^{-1}$ nor in $\omega$. Since $\gamma_i^{-1}(\mathbb{F}_2 - A_i) = B_i$ and also $\gamma_i^{-1} \bullet y = (\gamma_i^{-1}\zeta) \bullet x$, we conclude that $\gamma_i^{-1} \bullet y \in B_i^{\mathcal{O}}$, provided that $\gamma_i^{-1}\zeta$ does not end neither in $\tau^{-1}$ nor in $\omega$. Let us divide again into three subcases.

If $i = 0$, then $\gamma_0^{-1} = \sigma^{-n+1}$, which does not contain neither $\tau$ nor $\tau^{-1}$. Hence it is immediate to see that $\gamma_0^{-1}\zeta$ does not end neither in $\tau^{-1}$ nor in $\omega$ (which begins with $\tau$), as the same properties hold for $\zeta$.

If $i = 1$, then $\gamma_1^{-1} = \sigma^{-n+2}\tau^{-1}$. Since $\zeta$ does not end in $\omega$ and $\omega$ begins with $\tau$, then $\gamma_1^{-1}\zeta = \sigma^{-n+2}\tau^{-1}\zeta$ does not end in $\omega$. Moreover, if $\gamma_1^{-1}\zeta = \sigma^{-n+2}\tau^{-1}\zeta$ ends in $\tau^{-1}$, by the fact that $\zeta$ does not it follows that $\zeta = \varepsilon$. In this case $\zeta = \varepsilon \in A_1$ (by Lemma 2.6), contradicting the fact that $\zeta \in \mathbb{F}_2 - A_1$.

If $2 \leqslant i < n$, then $\gamma_i^{-1} = \sigma^{-i+2}\tau^{-1}\sigma^{i-1}$ and, as above, it is easy to see that $\gamma_i^{-1}\zeta$ does not end in $\omega$ if $\zeta$ does not. Finally, we have to exclude the case in which $\gamma_i^{-1}\zeta = \sigma^{-i+2}\tau^{-1}\sigma^{i-1}\zeta$ ends in $\tau^{-1}$. Since $\zeta$ does not end in $\tau^{-1}$, the only possibility is $\zeta = \sigma^{-i+1}$. In such case, $\zeta$ belongs to $A_i$ (by Lemma 2.6), contradicting the property that $\zeta \in \mathbb{F}_2 - A_i$. □

**Lemma 4.7.** *Let $\mathbb{F}_2 \curvearrowright \Omega$ locally abelian, $\mathcal{O}$ be an orbit with nontrivial stabilisers, and let $\omega \neq \varepsilon$ be of minimal length fixing at least one element $x \in \mathcal{O}$.*

*If $\omega = \sigma^{-k}$, for some $k > 0$, there is a partition*
$$\mathcal{O} = A_0^{\mathcal{O}} \sqcup \cdots \sqcup A_{n-1}^{\mathcal{O}} \sqcup B_0^{\mathcal{O}} \sqcup \cdots \sqcup B_{n-1}^{\mathcal{O}},$$
*with the property that $\gamma_i(B_i^{\mathcal{O}}) = \mathcal{O} - A_i^{\mathcal{O}}$ for every $i < n$.*

*Proof.* By the construction of Theorem 2.5 there are $A_i$ and $B_i$, for $0 \leqslant i < n$, with
$$\mathbb{F}_2 = A_0 \sqcup B_0 \sqcup \ldots \sqcup A_{n-1} \sqcup B_{n-1}$$
and $\varepsilon$ and $\omega$ both belong to the same element of the partition. Since $\omega = \sigma^{-k}$, the element of the partition to which they belong is $B_0$.

Now we define the sets $A_i^{\mathcal{O}}$ and $B_i^{\mathcal{O}}$ in a way that is analogous to Lemma 4.6: indeed, for $y \in \mathcal{O}$, let $\zeta$ be the unique word determined by Lemma 4.5, and we put $y \in A_i^{\mathcal{O}}$ or $y \in B_i^{\mathcal{O}}$ if and only if $\zeta \in A_i$ or $\zeta \in B_i$, respectively. By uniqueness of $\zeta$ and by the fact that the $A_i$ and $B_i$ form a partition of $\mathbb{F}_2$, we have
$$\mathcal{O} = A_0^{\mathcal{O}} \sqcup \ldots \sqcup A_{n-1}^{\mathcal{O}} \sqcup B_0^{\mathcal{O}} \sqcup \ldots \sqcup B_{n-1}^{\mathcal{O}}.$$
As above, it remains to prove that $\gamma_i(B_i^{\mathcal{O}}) = \mathcal{O} - A_i^{\mathcal{O}}$ for every $i < n$.

Let us prove $\gamma_i(B_i^{\mathcal{O}}) \subseteq \mathcal{O} - A_i^{\mathcal{O}}$ for $i < n$. Let $y \in B_i^{\mathcal{O}}$, then there is some $\zeta \in B_i$ that does not end neither in $\sigma$ nor in $\omega = \sigma^{-k}$ and such that $y = \zeta \bullet x$. Since
$$\gamma_i(B_i) = \mathbb{F}_2 - A_i,$$
we can write $\gamma_i \bullet y = (\gamma_i \zeta) \bullet x$, with $\gamma_i \zeta \in \mathbb{F}_2 - A_i$. In particular, if $\gamma_i \zeta$ does not end neither in $\sigma$ nor in $\sigma^{-k}$, then we conclude that $\gamma_i \bullet y \in \mathcal{O} - A_i^{\mathcal{O}}$.

If $i = 0$, then $\gamma_0 = \sigma^{n-1}$ and it is clear that $\gamma_0 \zeta = \sigma^{n-1} \zeta$ cannot end in $\sigma^{-k}$. However, it is possible that this word ends in $\sigma$. If this is the case, then we can consider $\sigma^{n-1} \zeta \omega = \sigma^{n-1} \zeta \sigma^{-k}$, which again is certainly not ending in $\sigma^{-k}$, since we are studying the case where $\sigma^{n-1}\zeta$ ends in $\sigma$. If this new word does not end in $\sigma$, then we stop. Otherwise, we repeat the argument until we find the first word that does not end in $\sigma$. Summing up, we find a word of the form $\sigma^{n-1}\zeta(\sigma^{-k})^a$, where $a \geqslant 0$, which does not end in $\sigma^{-k}$ nor in $\sigma$ and, since $\omega = \sigma^{-k}$ fixes $x$,

$$\gamma_0 \bullet y = (\sigma^{n-1}\zeta) \bullet x = (\sigma^{n-1}\zeta(\sigma^{-k})^a) \bullet x.$$

If we prove that $\sigma^{n-1}\zeta(\sigma^{-k})^a \in \mathbb{F}_2 - A_0$, then we get that $\gamma_0 \bullet y \in \mathcal{O} - A_0^{\mathcal{O}}$, reaching our goal. To prove that, notice that we are dealing with the case in which $\sigma^{n-1}\zeta$ ends in $\sigma$ (otherwise there is no need for this construction). Since $\zeta$ does not end in $\sigma$, then $\zeta$ is a (negative) power of $\sigma$. Hence $\sigma^{n-1}\zeta(\sigma^{-k})^a$ is itself a power of $\sigma$ and, by the construction in Theorem 2.5, it belongs to $B_0 \subseteq \mathbb{F}_2 - A_0$.

If $i = 1$, then $\gamma_1 = \tau\sigma^{n-2}$. Again it is easy to see that $\gamma_1 \zeta = \tau\sigma^{n-2}\zeta$ does not end in $\omega = \sigma^{-k}$ if $\zeta$ does not. Moreover since $\omega = \sigma^{-k}$ then, by the construction in Theorem 2.5, $B_1 = B_1^*$, so $\zeta \in B_1 = B_1^*$ begins with $\sigma^{-n+2}\tau^{-1}$. Therefore $\gamma_1 \zeta = \tau\sigma^{n-2}\zeta$ is a final part of $\zeta$: hence, it does not end in $\sigma$.

Finally we consider the case $2 \leqslant i < n$, where $\gamma_i = \sigma^{-i+1}\tau\sigma^{i-2}$. Again, by the construction in Theorem 2.5, we get $B_i = B_i^*$, so that $\zeta$ begins with $\sigma^{-i+2}\tau^{-1}$, meaning that $\zeta = \sigma^{-i+2}\tau^{-1}\eta$ for some $\eta \in \mathbb{F}_2$ not starting with $\tau$. Hence

$$\gamma_i \zeta = \sigma^{-i+1}\eta.$$

Since $\zeta$ does not end in $\sigma$, neither $\eta$ does. Hence also this word $\gamma_i \zeta$ cannot end in $\sigma$. If $\gamma_i \zeta$ does not end in $\sigma^{-k} = \omega$, we are done. Otherwise, we argue as in the case $i = 0$: we consider the word $\gamma_i \zeta (\sigma^k)^a = \sigma^{-i+1}\eta(\sigma^k)^a$, where $a \geqslant 1$ is the smallest such that this new word does not end in $\sigma^{-k}$. Also, by minimality of $a$, it does not end in $\sigma$. Since $\omega = \sigma^{-k}$ fixes $x$,

$$\gamma_i \bullet y = (\sigma^{-i+1}\eta) \bullet x = (\sigma^{-i+1}\eta(\sigma^k)^a) \bullet x$$

and so it suffices to prove that $(\sigma^{-i+1}\eta(\sigma^k)^a) \in \mathbb{F}_2 - A_i = \mathbb{F}_2 - A_i^*$. This is true since $\eta$ does not start with $\tau$, so that $(\sigma^{-i+1}\eta(\sigma^k)^a)$ does not begin with $\sigma^{-i+1}\tau$.

Then we are left to prove that $\gamma_i(B_i^{\mathcal{O}}) \supseteq \mathcal{O} - A_i^{\mathcal{O}}$ for $i < n$. This is equivalent to

$$\gamma_i^{-1}(\mathcal{O} - A_i^{\mathcal{O}}) \subseteq B_i^{\mathcal{O}}.$$

Let $y \in \mathcal{O} - A_i^{\mathcal{O}}$. This implies $y = \zeta \bullet x$ for some $\zeta \in \mathbb{F}_2 - A_i$ which does not end neither in $\sigma$ nor in $\sigma^{-k}$. Since $\gamma_i^{-1}(\mathbb{F}_2 - A_i) = B_i$ and also $\gamma_i^{-1} \bullet y = (\gamma_i^{-1}\zeta) \bullet x$, we deduce $\gamma_i^{-1} \bullet y \in B_i^{\mathcal{O}}$, provided that $\gamma_i^{-1}\zeta$ does not end neither in $\sigma$ nor in $\sigma^{-k}$.

If $i = 0$, then $\gamma_0^{-1} = \sigma^{-n+1}$ and it is clear that $\gamma_0^{-1}\zeta = \sigma^{-n+1}\zeta$ does not end in $\sigma$. If it does not end with $\sigma^{-k}$ either, then we are done. Otherwise, we can consider the word $\sigma^{-n+1}\zeta(\sigma^k)^a$, where $a \geqslant 1$ is the smallest natural number such that the corresponding word does not end in $\sigma^{-k}$. By minimality, such a word does not end in $\sigma$ either. Moreover, since $\omega = \sigma^{-k}$ fixes $x$, we can write

$$\gamma_0^{-1} \bullet y = (\sigma^{-n+1}\zeta) \bullet x = (\sigma^{-n+1}\zeta(\sigma^k)^a) \bullet x.$$

If we prove that $\sigma^{-n+1}\zeta(\sigma^k)^a \in B_0$, then we are done. Notice that $\sigma^{-n+1}\zeta$ ending in $\sigma^{-k}$ implies that $\zeta$ is a power of $\sigma$ (by the fact that $\zeta$ does not end neither in $\sigma$ nor $\sigma^{-k}$). Therefore the same is true of $\sigma^{-n+1}\zeta(\sigma^k)^a$. By construction, the latter belongs to $B_0$.

If $i = 1$, then $\gamma_1^{-1} = \sigma^{-n+2}\tau^{-1}$. By assumption $\zeta \in \mathbb{F}_2 - A_1 = \mathbb{F}_2 - A_1^*$ (by Theorem 2.5). This means that $\zeta$ does not begin with $\tau$: hence $\gamma_1^{-1}\zeta = \sigma^{-n+2}\tau^{-1}\zeta$ has $\zeta$ as final part and, like $\zeta$ itself, it does not end neither with $\sigma$ nor with $\sigma^{-k}$.

Finally, we deal with the case $2 \leqslant i < n$, where $\gamma_i^{-1} = \sigma^{-i+2}\tau^{-1}\sigma^{i-1}$. Again, by the construction in Theorem 2.5, we get $A_i = A_i^*$, so that $\zeta$ does not begin with $\sigma^{-i+1}\tau$. We notice that the word $\gamma_i^{-1}\zeta = \sigma^{-i+2}\tau^{-1}\sigma^{i-1}\zeta$ does not end in $\sigma^{-k}$: this follows from the fact that $\zeta$ does not end in $\sigma^{-k}$ and it does not begin with $\sigma^{-i+1}\tau$. If $\gamma_i^{-1}\zeta$ does not end in $\sigma$ either, then we conclude. Otherwise, we argue as in the case $i = 0$: we consider the word
$$\gamma_i^{-1}\zeta(\sigma^{-k})^a = \sigma^{-i+2}\tau^{-1}\sigma^{i-1}\zeta(\sigma^{-k})^a,$$
where $a \geqslant 1$ is the least natural number such that this new word does not end in $\sigma$. Also, by minimality of $a$, it does not end in $\sigma^{-k} = \omega$. By the fact that $\omega = \sigma^{-k}$ fixes $x$, we have
$$\gamma_i^{-1} \bullet y = (\sigma^{-i+2}\tau^{-1}\sigma^{i-1}\zeta) \bullet x = (\sigma^{-i+2}\tau^{-1}\sigma^{i-1}\zeta(\sigma^{-k})^a) \bullet x.$$
So it suffices to prove that $\sigma^{-i+2}\tau^{-1}\sigma^{i-1}\zeta(\sigma^{-k})^a \in B_i = B_i^*$. This is true since $\zeta$ does not start with $\sigma^{-i+1}\tau$, so $\sigma^{-i+2}\tau^{-1}\sigma^{i-1}\zeta(\sigma^{-k})^a$ begins with $\sigma^{-i+2}\tau^{-1}$. □

Hence we can now prove the following result:

**Corollary 4.8** (Generalised Banach-Tarski for $\mathbb{S}^2$). *With respect to $\mathrm{SO}(3) \curvearrowright \mathbb{R}^3$, the sphere $\mathbb{S}^2$ is $(n, 2n)$-paradoxical and not $(n, r)$-paradoxical for $r < 2n$.*

*Proof.* By Fact 4.1, the group $\mathrm{SO}(3)$ has a subgroup isomorphic to $\mathbb{F}_2$. By Corollary 3.8, $\mathrm{SO}(3)$ is $(n, 2n)$-paradoxical. The action $\mathrm{SO}(3) \curvearrowright \mathbb{S}^2$ is locally abelian because of Lemma 4.3, hence by Theorem 4.4 we get the $(n, 2n)$-paradoxality of $\mathbb{S}^2$. Moreover, Theorem 3.4 implies that $\mathbb{S}^2$ cannot be $(n, r)$-paradoxical for $r < 2n$. □

## 5. The Result for the Ball

In this section we apply our generalisation of the Banach-Tarski result for the sphere $\mathbb{S}^2$ (Corollary 4.8) to find the minimal number of pieces needed to perform the same construction with the closed ball $\mathbb{D}^3$.

**Definition 5.1.** A *closed hemisphere* of $\mathbb{S}^2$ is a subset $\mathcal{H} \subseteq \mathbb{S}^2$ for which there exist $a, b, c \in \mathbb{R}$ not all zero such that
$$\mathcal{H} = \{\langle x, y, z\rangle \in \mathbb{S}^2 : ax + by + cz \leqslant 0\}.$$
If in the definition of $\mathcal{H}$ we replace $\leqslant$ with $<$, we say that $\mathcal{H}$ is an *open hemisphere*. Clearly, one can define closed and open hemispheres relative to any sphere in $\mathbb{R}^3$ analogously.

**Lemma 5.2.** *Let $\mathcal{B} \subseteq \mathbb{R}^3$ be a closed ball of radius $1$ with center $\langle a, b, c\rangle \neq \langle 0, 0, 0\rangle$. Then there is a closed hemisphere of $\mathbb{S}^2$ which is disjoint from $\mathcal{B}$.*

*Proof.* Consider the closed hemisphere defined by
$$\mathcal{H} := \{\langle x, y, z\rangle \in \mathbb{S}^2 : ax + by + cz \leqslant 0\}.$$
Let $\langle x, y, z\rangle \in \mathcal{H} \cap \mathcal{B}$. By definition of $\mathcal{B}$, we have
$$(x - a)^2 + (y - b)^2 + (z - c)^2 \leqslant 1.$$

Since $x^2 + y^2 + z^2 = 1$ we get
$$a^2 + b^2 + c^2 \leqslant 2(ax + by + cz) \leqslant 0,$$
which is a contradiction because $\langle a, b, c \rangle \neq \langle 0, 0, 0 \rangle$, hence $a^2 + b^2 + c^2 > 0$. $\square$

**Definition 5.3.** A map $\sigma : \mathbb{R}^3 \to \mathbb{R}^3$ is called a *direct isometry* if it is of the form
$$\sigma(\mathbf{x}) = \rho(\mathbf{x}) + \mathbf{b}$$
for some $\rho \in \mathrm{SO}(3)$ and some $\mathbf{b} \in \mathbb{R}^3$.

We denote by $\mathrm{SE}(3)$ the group of direct isometries on $\mathbb{R}^3$.

**Theorem 5.4.** *Assume that $\mathbb{D}^3$ is $(n, k)$-paradoxical with respect to $\mathrm{SE}(3) \curvearrowright \mathbb{R}^3$. Then $k \geqslant 3n - 1$.*

*Proof.* By assumption, we have $\mathbb{D}^3 = P_0 \sqcup \ldots \sqcup P_{n-1}$ with $P_j \sim_{k_j} \mathbb{D}^3$ for all $j < n$ and $k_0 + \cdots + k_{n-1} = k$. Assume by contradiction $k < 3n - 1$. Since it is not possible that $P_j \sim_1 \mathbb{D}^3$, there exist $0 \leqslant j_0 < j_1 < n$ such that $k_{j_0} = k_{j_1} = 2$.

Hence there are pairwise disjoint subsets $B_1$, $B_2$, $B_3$, $B_4$ of $\mathbb{D}^3$ and isometries $\sigma_1, \sigma_2, \sigma_3, \sigma_4$ with
$$\sigma_1[B_1] \sqcup \sigma_2[B_2] = \sigma_3[B_3] \sqcup \sigma_4[B_4] = \mathbb{D}^3,$$
where, here and in the following, $\sigma[X]$ denotes the image of the set $X$ under $\sigma$.

We prove that at least one of the four maps $\sigma_1$, $\sigma_2$, $\sigma_3$ and $\sigma_4$ does not fix the origin $\mathbf{0} := \langle 0, 0, 0 \rangle$ of $\mathbb{R}^3$. Indeed, assume this is not the case. By assumption, at most one among $B_1$, $B_2$, $B_3$ and $B_4$ contains $\mathbf{0}$. Without loss of generality, we can assume $\mathbf{0} \notin B_3 \cup B_4$. Then, since $\sigma_3(\mathbf{0}) = \mathbf{0}$ and $\sigma_4(\mathbf{0}) = \mathbf{0}$, we get
$$\mathbf{0} \notin \sigma_3[B_3] \cup \sigma_4[B_4] = \mathbb{D}^3.$$

Consequently, this proves that there is some $1 \leqslant i \leqslant 4$ such that $\sigma_i(\mathbf{0}) \neq \mathbf{0}$. Without loss of generality assume that $i = 4$. Then $\sigma_4[\mathbb{D}^3]$ is a closed ball with center different from $\mathbf{0}$. By Lemma 5.2, there is a closed hemisphere $\mathcal{H}$ of $\mathbb{S}^2$ disjoint from $\sigma_4[\mathbb{D}^3]$. Since $\sigma_3[B_3] \cup \sigma_4[B_4] = \mathbb{D}^3$, then $\mathcal{H} \subseteq \sigma_3[B_3]$, so $\sigma_3^{-1}[\mathcal{H}] \subseteq B_3$.

Now, since $\mathcal{H}$ is a closed hemisphere of $\mathbb{S}^2$, then $\sigma_3^{-1}[\mathcal{H}]$ is a closed hemisphere of some sphere of radius 1. Since
$$\sigma_3^{-1}[\mathcal{H}] \subseteq B_3 \subseteq \mathbb{D}^3,$$
then $\sigma_3^{-1}[\mathcal{H}]$ is a closed hemisphere of $\mathbb{S}^2$.

Since $B_1$ and $B_2$ are disjoint from $B_3$, it follows that $(B_1 \cup B_2) \cap \mathbb{S}^2$ is contained in the complement of $\sigma_3^{-1}[\mathcal{H}]$ in $\mathbb{S}^2$, which is an open hemisphere of $\mathbb{S}^2$. In particular, neither $B_1$ nor $B_2$ contain a closed hemisphere of $\mathbb{S}^2$.

If $\sigma_1(\mathbf{0}) \neq \mathbf{0}$, by the same argument that we used for $B_3$, we get that $B_2$ contains a closed hemisphere of $\mathbb{S}^2$. Hence $\sigma_1(\mathbf{0}) = \mathbf{0}$ and the same argument applied to $\sigma_2$ shows that $\sigma_2(\mathbf{0}) = \mathbf{0}$. Therefore $\sigma_1$ and $\sigma_2$ are direct isometries fixing $\mathbf{0}$. Thus they are rotations and, in particular, they map $\mathbb{S}^2$ into itself. Hence we have
$$(\sigma_1[B_1] \cup \sigma_2[B_2]) \cap \mathbb{S}^2 = \sigma_1[B_1 \cap \mathbb{S}^2] \cup \sigma_2[B_2 \cap \mathbb{S}^2].$$

As we have already pointed out, $B_1 \cap \mathbb{S}^2$ and $B_2 \cap \mathbb{S}^2$ are contained in an open hemisphere of $\mathbb{S}^2$, hence the same holds for $\sigma_1[B_1 \cap \mathbb{S}^2]$ and $\sigma_2[B_2 \cap \mathbb{S}^2]$. Therefore the equality above implies that $(\sigma_1[B_1] \cup \sigma_2[B_2]) \cap \mathbb{S}^2$ is contained in the union of two open hemispheres, so it is strictly contained in $\mathbb{S}^2$: this is a contradiction with the assumption that $\sigma_1[B_1] \sqcup \sigma_2[B_2] = \mathbb{D}^3$

Hence the initial assumption $k < 3n - 1$ must be rejected. □

**Theorem 5.5.** $\mathbb{D}^3$ *is* $(n, 3n-1)$*-paradoxical with respect to* $\mathrm{SE}(3) \curvearrowright \mathbb{R}^3$.

*Proof.* For $0 < r \leqslant 1$, let $S_r := \{\mathbf{x} \in \mathbb{R}^3 : \|\mathbf{x}\| = r\}$ be the sphere of radius $r$ in $\mathbb{R}^3$ centered at $\mathbf{0}$. Notice that $\mathbb{D}^3 = \{\mathbf{0}\} \cup \bigcup_{0 < r \leqslant 1} S_r$ and $S_1 = \mathbb{S}^2$.

Recall that $\mathrm{SO}(3)$ contains a subgroup $F$ isomorphic to the free group with two generators $\mathbb{F}_2$. We identify $F$ with $\mathbb{F}_2$ and we denote by $\gamma_i$, $0 \leqslant i < n$, the rotations corresponding to the reduced words defined in Section §2. Similarly, we denote by $\varepsilon$, $\sigma^{-1}$, ..., $\sigma^{-n+2}$ the rotations associated to the words with the same names by a fixed isomorphism between $F$ and $\mathbb{F}_2$.

Next, recall that the action of $F = \mathbb{F}_2$ on $\mathbb{S}^2$ is locally abelian. Hence, for each $0 < r < 1$, by Corollary 4.8 (in particular, by the construction in the proof of Theorem 4.4), we can split $S_r$ in $2n$ pieces, which we call $A_0^r, B_0^r, \ldots, A_{n-1}^r, B_{n-1}^r$, such that $\gamma_i(B_i^r) = S_r - A_i^r$ for all $0 \leqslant i < n$. Here we point out that Corollary 4.8 holds for every sphere, by homothety.

We assume that there exists a partition of $S_1$ into $3n-1$ pieces, that we denote by $A_0^1, B_0^1, \ldots, A_{n-1}^1, B_{n-1}^1$ and $\{\mathbf{x}_1\}, \{\mathbf{x}_2\}, \ldots, \{\mathbf{x}_{n-1}\}$, such that $\gamma_i(B_i^1) = S_1 - A_i^1$, for all $0 \leqslant i < n$, and that each $\{\mathbf{x}_i\}$ is a singleton. We postpone the proof of existence of such a partition and we complete the proof of the theorem first.

To such an extent, we consider a partition of $\mathbb{D}^3$ consisting of $3n - 1$ pieces

$$\{A_0, \{B_i\}_{0 \leqslant i < n}, \{A_i\}_{1 \leqslant i < n}, \{\mathbf{x}_i\}_{1 \leqslant i < n}\}$$

where we define $A_0$, $B_i$ for $0 \leqslant i < n$ and $A_i$ for $1 \leqslant i < n$ as:

$$A_0 := \{\mathbf{0}\} \cup \bigcup_{0 < r \leqslant 1} A_0^r; \quad B_i := \bigcup_{0 < r \leqslant 1} B_i^r; \quad A_i := \bigcup_{0 < r \leqslant 1} A_i^r.$$

Let $\varrho_i$ be the Euclidean translation which maps $\mathbf{x}_i$ to $\mathbf{0}$. Then we have

$$\begin{aligned} A_0 \cup \gamma_0(B_0) &= A_0 \cup \gamma_0(\textstyle\bigcup_{0 < r \leqslant 1} B_0^r) \\ &= A_0 \cup \textstyle\bigcup_{0 < r \leqslant 1} \gamma_0(B_0^r) \\ &= \{\mathbf{0}\} \cup \textstyle\bigcup_{0 < r \leqslant 1} A_0^r \cup \textstyle\bigcup_{0 < r \leqslant 1} (S_r - A_0^r) \\ &= \{\mathbf{0}\} \cup \textstyle\bigcup_{0 < r \leqslant 1} S_r = \mathbb{D}^3. \end{aligned}$$

Similarly, for $0 < i < n$ we have that

$$\begin{aligned} A_i \cup \gamma_i(B_i) \cup \varrho_i(\{\mathbf{x}_i\}) &= A_i \cup \gamma_i(\textstyle\bigcup_{0 < r \leqslant 1} B_i^r) \cup \{\mathbf{0}\} \\ &= \textstyle\bigcup_{0 < r \leqslant 1} A_i^r \cup \textstyle\bigcup_{0 < r \leqslant 1} \gamma_i(B_i^r) \cup \{\mathbf{0}\} \\ &= \textstyle\bigcup_{0 < r \leqslant 1} (A_i^r \cup (S_r - A_i^r)) \cup \{\mathbf{0}\} = \mathbb{D}^3. \end{aligned}$$

This proves that $\mathbb{D}^3$ is $(n, 3n-1)$-paradoxical with respect to the action of $\mathrm{SE}(3)$ on $\mathbb{R}^3$. In particular, we have obtained $n-1$ copies of $\mathbb{D}^3$, each formed by 3 pieces, and another copy (the one constructed from $A_0$ and $B_0$), formed by 2 pieces.

Therefore, it remains to prove that the partition of $S_1 = \mathbb{S}^2$ with the above mentioned properties does exist. First of all, we consider the action of $F \simeq \mathbb{F}_2$ on $\mathbb{S}^2$ and we fix an orbit $\mathcal{O}$ with trivial stabilisers. Such an orbit does exist since $F$ is countable and each rotation different from the identity fixes only two points on $\mathbb{S}^2$.

Therefore there exists a point (hence an orbit) with trivial stabiliser. We partition the set $\mathbb{S}^2 - \mathcal{O}$ as we did for $S_r$, for $0 < r < 1$, by using Corollary 4.8.

On the orbit $\mathcal{O}$, we pick an element $\mathbf{x}_1$ and we define
$$\mathbf{x}_j := \sigma^{-j+1}(\mathbf{x}_1),$$
for $2 \leqslant j < n$. Since $\mathcal{O}$ has trivial stabiliser, for every $\mathbf{y} \in \mathcal{O} - \{\mathbf{x}_j : 1 \leqslant j < n\}$ there is a unique element
$$\omega \in F - \{\varepsilon, \sigma^{-1}, \ldots, \sigma^{-n+2}\}$$
such that $\mathbf{y} = \omega(\mathbf{x}_1)$. Hence we partition $\mathcal{O}$ in $2n$ subsets as we previously did: we let $\mathbf{y} \in A_i^1$ if $\omega \in A_i^*$ and $\mathbf{y} \in B_i^1$ if $\omega \in B_i^*$ for every $0 \leqslant i < n$, where $A_i^*, B_i^*$ are defined in Section §2. Now, by applying Lemma 2.3, we get the existence of a partition with the required properties. □

With reference to the construction of the sets $A_i^*$ and $B_i^*$ preceding Lemma 2.3, in addition to $\varepsilon$ we left behind $n - 2$ elements: those having the form $\sigma^{-i}$, for $0 < i < n - 1$. This is a difference with the classical result for $n = 2$, where in the corresponding partition every element of $\mathbb{F}_2$ except $\varepsilon$ is included. Our modification plays a crucial role in the proof of Theorem 5.5 in the construction of the $n - 1$ singleton pieces $\{\mathbf{x}_1\}, \ldots, \{\mathbf{x}_{n-1}\}$ of the desired partition.

**Corollary 5.6** (Generalised Banach-Tarski for $\mathbb{D}^3$). *With respect to* $\mathrm{SE}(3) \curvearrowright \mathbb{R}^3$, *the ball $\mathbb{D}^3$ is $(n, 3n - 1)$-paradoxical and is not $(n, r)$-paradoxical for $r < 3n - 1$.*

*Proof.* Combine Theorem 5.4 and Theorem 5.5. □

We close with the observation that the proof of Theorem 5.4 not only shows that $3n - 1$ pieces are necessary, but also that, in this minimal case, one of the $n$ "new" balls must be formed by two pieces and all of the others by three. The partition that we define in the proof of Theorem 5.5 satisfies such requirement: $n - 1$ of the pieces are singletons, each contributing to one of the balls formed by three pieces. Hence a more meaningful way to express the minimal number may be $2 + 3(n - 1)$.

## Acknowledgments

The authors are members of the Istituto Nazionale di Alta Matematica (INdAM): the first author is in the Gruppo Nazionale per le Strutture Algebriche, Geometriche e le loro Applicazioni (GNSAGA); the second author is in the Gruppo Nazionale per l'Analisi Matematica, la Probabilità e le loro Applicazioni (GNAMPA).

We both would like to thank deeply Stefano Baratella, who had also been our Bachelor's degrees supervisor, for his constant support and his careful and thorough reading of an earlier draft of this paper. His constructive comments helped improving the quality of the article.

Departament de Matemàtiques i Informàtica, Universitat de Barcelona, Gran Via de les Corts Catalanes, 585, 08007 Barcelona, Catalunya; Dipartimento di Matematica, Università degli Studi di Trento, Via Sommarive, 14, 38123 Povo (Trento), Italia
*Email address*: straffelini@ub.edu; cesare.straffelini@unitn.it

Dipartimento di Matematica, Università degli Studi di Trento, Via Sommarive, 14, 38123 Povo (Trento), Italia
*Email address*: kilian.zambanini@unitn.it